\documentclass[12pt]{amsart}%
\usepackage{amsfonts}
\usepackage{amsmath}
\usepackage{amssymb}
\usepackage{graphicx}
\usepackage{setspace}
\makeatletter
\@addtoreset{equation}{section}
\makeatother

\newtheorem{theorem}{Theorem}[section]

\newtheorem{corollary}[theorem]{Corollary}

\newtheorem{lemma}[theorem]{Lemma}

\newtheorem{remark}[theorem]{Remark}

\usepackage{xcolor}

\makeatletter
\def\@makefnmark{}
\makeatother

\begin{document}
\title[critical anisotropic Sobolev equation without the finite volume constraints]{ classification of positive solutions of critical anisotropic Sobolev equation without the finite volume constraint}

\author{Lu Chen}
\address{Key Laboratory of Algebraic Lie Theory and Analysis of Ministry of Education, School of Mathematics and Statistics, Beijing Institute of Technology, Beijing, 100081, P. R. China}
\email{chenlu5818804@163.com}

\author{Tian Wu}
\address{School of Mathematical Sciences, University of Science and Technology of China, Hefei, Anhui, 230026, P. R. China}
\email{wt1997@ustc.edu.cn}

\author{Jin Yan}
\address{School of Mathematical Sciences, University of Science and Technology of China, Hefei, Anhui, 230026, P. R. China}
\email{yjoracle@mail.ustc.edu.cn}

\author{Yabo Yang}
\address{Key Laboratory of Algebraic Lie Theory and Analysis of Ministry of Education, School of Mathematics and Statistics, Beijing Institute of Technology, Beijing, 100081, P. R. China}
\email{Yabo\underline{~}Yang0927@outlook.com}

\thanks{The first author was partly supported by the National Key Research and Development Program (No.
2022YFA1006900) and National Natural Science Foundation of China (No. 12271027),
the second author was supported by Anhui Postdoctoral Scientific Research Program Foundation (Grant No. 2025B1055) and the Fundamental Research Funds for the Central Universities (Grant No. WK0010250106), the third author was supported by Anhui Postdoctoral Scientific Research Program Foundation (Grant No. 2025B1055), the fourth author was supported by a grant from Beijing Institute of Technology (No. 2022CX01002).}

\begin{abstract}
In this paper, we classify all positive solutions of the critical anisotropic Sobolev equation
\begin{equation}\label{0.1}
-\Delta^{H}_{p}u = u^{p^{*}-1}, \ \ x\in \mathbb{R}^n
\end{equation}
without the finite volume constraint for $n \geq 3$ and $p_n(\Lambda) < p < n$, where $p^{*} = \frac{np}{n-p}$ denotes the critical Sobolev exponent, $-\Delta^{H}_{p}=-div(H^{p-1}(\cdot)\nabla H(\cdot))$ denotes the anisotropic $p$-Laplace operator and $\Lambda = \lambda\max\limits_{\substack{\xi \in \mathbb{R}^n\\1 \leq i, j \leq n}}\left\{\frac{|\xi|^{2}(\nabla^{2}_{ij}H^{p}(\xi))} {p(p-1)H^{p}(\xi)}\right\}$. By employing a novel approach based on invariant tensors technique, and using a Kato-type inequality, we prove that the positive solutions of \eqref{0.1} can be classified for $p_n(\Lambda) \leq p < n$, where $p_n(\Lambda)$ depends explicitly on $\Lambda$. This result removes the finite volume assumption on the classification of critical anisotropic $p$-Laplace equation which was obtained by Ciraolo-Figalli-Roncoroni in the literature \cite{CFR}. In particular, this results capture the precise dependence of critical exponents $p$ on both $n$ and $\Lambda$.
\end{abstract}

\maketitle {\small {\bf Keywords:} Critical anisotropic Sobolev equation, Classification, Without finite volume constraint, Integral inequality, Regularity, Invariant tensor technique\\

}
\section{Introduction}
Given $n \geq 2$ and $1<p<n$, the classical Sobolev inequality \cite{Sobolev} in $\mathbb{R}^n$ states that for any $u\in W^{1,p}(\mathbb{R}^n)$, there holds
\begin{equation}\label{1.1}
\int_{\mathbb{R}^n}|u|^{p^{*}}dx \leq C(n, p, s)\int_{\mathbb{R}^n}|\nabla u |^{p}dx,
\end{equation}
where $p^{*} = \frac{np}{n-p}$ denotes the critical Sobolev exponent. Aubin \cite{Aubin} and Talenti \cite{Talenti} applied the technique of symmetry and rearrangement combining the Bliss Lemma to show that all radial extremals of Sobolev inequality must take the form as
\begin{equation*}
U=\big(1+|x|^{\frac{p}{p-1}}\big)^{-\frac{n-p}{p}},
\end{equation*}
up to some dilation and translation. However, they didn't classify  all extremals of Sobolev inequality. Later, Erausquin, Nazaret and Villani \cite{ENV} showed that all extremals must take the form as
\begin{equation*}
U=\big(1+|x|^{\frac{p}{p-1}}\big)^{-\frac{n-p}{p}},
\end{equation*}
up to some dilation and translation by the optimal transportation method. Obviously, the extremals of Sobolev inequality satisfy the critical Sobolev equation:
\begin{equation}\label{1.2}
\begin{cases}
-\Delta_{p}u = u^{\frac{np}{n-p}-1} &x \in \mathbb{R}^n,\\
u\geq 0 &x \in \mathbb{R}^n,\\
u \in\;W^{1,p}(\mathbb{R}^n).
\end{cases}
\end{equation}
The classification of positive solutions of equation \eqref{1.2} started in the crucial papers \cite{GNN} and \cite{GS} and it has been the object of several studies. Damascelli-Merch\'{a}n-Montoro-Sciunzi \cite{DMMS}, Sciunzi \cite{Sciunzi} and V\'{e}tois \cite{Vet16} established the symmetry of positive solutions of equation \eqref{1.2}, which together with Aubin and Talenti's results deduces the uniqueness of extremals of Sobolev inequality. And Jerison-Lee \cite{JL} employed computer-assisted calculations to prove a classification theorem with the assumption of finite energy.
\vskip0.1cm
A natural problem is whether we can classify the positive solutions of critical Sobolev equations \eqref{1.2} without the finite volume assumption. In fact, this is proved to be true by Caffralli-Giddas-Spruck \cite{CGS} applying Kelvin transform and moving plane method to classify all the positive solutions of the Yamabe equation \cite{LP} when $p=2$. Later, Chen-Li \cite{CL} provided a simpler proof using the moving plane method to obtain the same results. We also note that Dai-Liu-Qin \cite{DLQ} and Dai-Qin \cite{DQ18, DQ21} applied the method of moving spheres in integral form to classify all nonnegative solutions to the integral equations, the conformally invariant system with mixed order and exponentially increasing nonlinearity and the high-order equations, respectively. Later, Peng \cite{Peng} applying the same method to classify the solutions to mixed order elliptic system with general nonlinearity. Beyond positive solutions, we also mention that some classification results regarding sign-changing solutions to the equation
\begin{align*}
-\Delta_{p}u = |u|^{\alpha-1}u \;\;x \in \mathbb{R}^n,
\end{align*}
have been classified by Bahri-Lions \cite{BL} for $p = 2$. Subsequently, Farina \cite{Farina} and Damascelli-Farina-Sciunzi-Valdinoci \cite{DFSV} classified stable solutions for $p = 2$ and $p >2$. Furthermore, Farina-Sciunzi-Vuono \cite{FSV} studied the established corresponding Liouville theorems for stable solutions to the more general quasilinear equation.
\par However, the Kelvin transform is not available for the general $p$-Laplace equation, hence the classification problem of critical Sobolev equation for $p\neq 2$ without the finite volume assumption is a challenging problem. Recently, Catino-Monticelli-Roncoroni \cite{CMR} solved the classification problem under the assumption $\frac{n}{2} < p <2$ in $n=2,3$, Ou \cite{Ou} for $\frac{n+1}{3} \leq p<n$ and V\'{e}tois \cite{Vet24} for $p_n < p <n$ where
\begin{equation}\label{1.3}
p_n=
\begin{cases}
\frac{8}{5}&\mathrm{if}\;\;\;n=4,\\
\frac{4n+3-\sqrt{4n^2+12n-15}}{6}&\mathrm{if}\;\;\;n\geq 5 .
\end{cases}
\end{equation}
 The same method has been also used successfully in the analogous problems such as critical Sobolev on Euclidean space, Heisenberg group and \emph{C-R} manifold (see \cite{CMR}, \cite{MO}, \cite{MOW} and \cite{Vet24}). However, the classification result for critical Sobolev equation without the finite volume assumption in the remaining index still keeps open.
\vskip0.1cm
Now, let us turn to the introduction of anisotropic Sobolev inequality. Anisotropic Sobolev inequality can be stated as follows: for any $u\in W^{1,p}(\mathbb{R}^n)$, there holds
\begin{equation}\label{1.4}
\left(\int_{\mathbb{R}^n}|u|^{p^{*}}dx\right)^{\frac{1}{p^{*}}} \leq C(n, p)\left(\int_{\mathbb{R}^n}H^{p}(\nabla v)dx\right)^{\frac{1}{p}},
\end{equation}
where $H$ is a 1-homogenous convex function ( see details in subsection 2.1) and $C_{n,p}$ denotes the best possible constant which makes the anisotropic Sobolev inequality holds. This sharp inequality was first obtained by Alvino-Ferone-Trombetti-Lions \cite{AFTL} using the convex symmetrization technique.  However, they did not solve the uniqueness problem of extremals of anisotropic Sobolev inequality. Ciraolo, Figalli and Roncoroni [\cite{CFR}, Appendix A] solved the uniqueness problem by adapting the optimal transportation method. Furthermore, they proved that all positive solutions of anisotropic Sobolev equation with the finite volume constraint
\begin{equation}\label{1.5}
\begin{cases}
-\mathrm{div}\;(a(\nabla u)) = u^{p^{*}-1}&x \in \mathbb{R}^n,\\
u\geq 0&x \in \mathbb{R}^n,\\
\int_{\mathbb{R}^n}|u|^{\frac{np}{n-p}}dx<+\infty,
\end{cases}
\end{equation}
must take the form as
\begin{equation*}
U_{\lambda}(x) = \left(\frac{(\lambda^{\frac{1}{p-1}}(n^{\frac{1}{p}}(\frac{n-p}{p-1})^{\frac{p-1}{p}})}{\lambda^{\frac{p}{p-1}}+H_0(x)^{\frac{p}{p-1}}}\right)^{\frac{n-p}{p}},
\end{equation*}
where $a(\nabla u) = H^{p-1}(\nabla u)\nabla H(\nabla u)$, up to some translation. They classify all positive solutions and furthermore extended the classification results to the case of critical anisotropic Sobolev equation in convex cone. Recently, Montoro-Muglia-Sciunzi \cite{MMS} classify all weak solutions to Laplacian equation in half space using the similar method.
\par
It should be noted that in the research of anisotropic Sobolev equation, the finite volume assumption plays an important role. In this paper, we are devoted to classify positive solutions of critical anisotropic Sobolev equation without the finite volume constraint:
\begin{equation}\label{1.6}\begin{cases}
-\mathrm{div}\;(a(\nabla u)) = u^{p^{*}-1}& x \in \mathbb{R}^n,\\
u\geq 0& x \in \mathbb{R}^n.
\end{cases}\end{equation}
We are motivated by recent progress in  Liang-Wu-Yan's work in \cite{LWY}, Ma-Ou-Wu's work \cite{MOW}, Ou's work in \cite{Ou} and V\'{e}tois's work in \cite{Vet24}. We have found that invariant tensor method in literature \cite{LWY} simplifies the computational process and for this reason we provide a proof by suitably adapting the invariant tensor method to classify the positive solutions of critical anisotropic Sobolev equation for the case of $p_n(\Lambda) < p < n$, which could capture the precise dependence of critical exponents $p$ on both $n$ and $\Lambda$. Our main result states as:
\medskip
\vskip0.1cm

\begin{theorem}\label{thm 1.1}
For $p_n(\Lambda) < p < n$ and $\Lambda < 1 + c(n)$, assume that $u \in W^{1,p}_{loc}(\mathbb{R}^n)$ is a positive weak solution of \eqref{1.6}. Then $u$ must take the form as
\begin{equation*}
U_{\lambda}(x) = \left(\frac{(\lambda^{\frac{1}{p-1}}(n^{\frac{1}{p}}(\frac{n-p}{p-1})^{\frac{p-1}{p}})}{\lambda^{\frac{p}{p-1}}+H_0(x)^{\frac{p}{p-1}}}\right)^{\frac{n-p}{p}},
\end{equation*}
up to some translation, where
\begin{equation}\label{p_n}
p_n=
\begin{cases}
\frac{8}{5}&\mathrm{if}\;\;\;n=4,\\
\frac{(5+n+3n\Lambda-2\Lambda)-\sqrt{(2\Lambda-n-5-3n\Lambda)^2-12(n^2\Lambda+n+2)}}{6}&\mathrm{if}\;\;\;n\geq 5,
\end{cases}
\end{equation}
and $c(n)$ is a constant depending on $n$ which could be expressed precisely in the proof of Theorem \ref{thm 1.1}.
\end{theorem}

\begin{remark}
The proof of Theorem \ref{thm 1.1} need to construct the vital integral inequality \eqref{3.3.2} involving the suitable vector $a^{i}$ and matrix $W$. Applying this integral inequality, the decay estimate in Lemma \ref{lem 3.8}, and Kato's type inequality in Lemma \ref{lem 3.6}, through complicated calculation, we could obtain that $Tr(W^2)$ is equal to zero. This together with construction of $W$ deduces $W=0$, which can help to classify all positive solutions of critical anisotropic Sobolev equation. It should be noted that when $\Lambda =1$, our results coincide with those obtained by V\'{e}tois in \cite{Vet24}, demonstrating the consistency of our approach with prior work, i.e. $p_n$ is equal to \eqref{1.3}. Moreover, our method extends their formulation and achieves improved performance in more general cases.
\end{remark}
\vskip0.1cm

\begin{remark}
A function $u \in W^{1,p}_{loc}(\mathbb{R}^n) \bigcap  L^{\infty}_{loc}(\mathbb{R}^n)$ is said to be a weak solution of \eqref{1.6} if
\begin{align}\label{1.7}
\int_{\mathbb{R}^n}H^{p-1}(\nabla u)\nabla H(\nabla u)\cdot\nabla \psi dx -\int_{\mathbb{R}^n}u^{p^{*}-1} \psi = 0,
\end{align}
for any $\psi \in C^{\infty}_{c}(\mathbb{R}^n)$.
\end{remark}
\vskip0.1cm
Here we mention some well-known facts about solutions of \eqref{1.6}, for any positive weak solution $u$ of \eqref{1.6}, we have
\begin{equation}\label{1.8}
u \geq C(n,p,\min_{|x|=1}u)|x|^{-\frac{n-p}{p-1}}    \;\;\;\;\;    \text{for} \;|x| > 1,
\end{equation}
where  $C$ is denote as a general positive constant. In fact, the estimate \eqref{1.8} has been derived for positive weak super $p$-harmonic functions (see \cite{CFR}).

\medskip
\vskip0.1cm
\noindent\textbf{Organization of the paper:}
This paper is organized as follows. In section 2 we introduce some notations involving anisotropic norms and provide a brief proof of regularity of solutions of critical anisotropic Sobolev equation. In Section 3, we construct suitable vector fields and establish the vital integral inequality \eqref{3.3.1} which plays a crucial role on classification of critical anisotropic Sobolev equation. In section 5 we provide a new approach by suitably adapting the invariant tensor method to classify the positive solutions of critical anisotropic Sobolev equation for $p_n(\Lambda) < p < n$, thereby capturing the precise dependence of critical exponents $p$ on both $n$ and $\Lambda$.

\section{Preliminaries}
\par In this section, we introduce some basic notations and properties about anisotropic norms and present the regularity of weak solutions of anisotropic equation. For more properties of anisotropic operators, we refer the readers to references \cite{BC}, \cite{CFR}, \cite{CY24}, \cite{CY25}, \cite{WX11}, \cite{WX12} and its references therein.
\medskip
\vskip0.1cm
\noindent\textbf{2.1 Some basic properties of anisotropic norms:}
Let $H: \mathbb{R}^n \rightarrow \mathbb{R}$ be a norm such that $H^2$ is of class $C^2(\mathbb{R}\backslash \{0\})$ and it is uniformly convex. This fact is easily seen to be equivalent to the following three properties:
\vskip0.1cm
\begin{align*}
H \;\;\mathrm{is} \;\;\mathrm{convex};
\end{align*}\[\]
\begin{align*}
H(\xi) \geq 0 \;\;\mathrm{for}\;\; \xi \in \mathbb{R}^n \;\;\mathrm{and}\;\;H(\xi) = 0 \;\;\mathrm{if}\;\mathrm{and}\;\mathrm{only}\;\mathrm{if}\;\;\xi = 0;
\end{align*}\[\]
\begin{align}\label{2.1.1}
H(t\xi) = tH(\xi)\;\;\mathrm{for}\;\; \xi \in \mathbb{R}^n \;\;\mathrm{and}\;\;\mathrm{for}\;\;t > 0.
\end{align}\[\]
All norms in $\mathbb{R}^n$ are equivalent. Hence, there exist positive constants $\lambda_1$ and $\lambda_2$ depending on $n, p, H$ such that
\vskip0.1cm
\begin{align}\label{2.1.2}
\lambda_1|\xi|^{p-2}|\zeta|^2 \leq \frac{1}{p}\nabla^2_{\xi_i\xi_j}H^{p}(\xi)\zeta_i\zeta_j \leq \lambda_2|\xi|^{p-2}|\zeta|^2 \;\;\mathrm{for}\;\;\xi\in \mathbb{R}^n.
\end{align}\[\]
Accordingly, $H_0$ denotes the dual norm to $H$ given by
\vskip0.1cm
\begin{align}\label{2.1.3}
H_0(\xi) = \sup\limits_{\xi\neq 0}\frac{\xi\cdot\eta}{H(\xi)} \;\;\;\forall \eta \in \mathbb{R}^n.
\end{align}\[\]
The following properties
\begin{align}\label{2.1.4}
H(\nabla_{\eta} H_0(\eta)) = 1,\;\;\;\;\;H_0(\nabla_{\xi}H(\xi)) = 1,\ \ \forall\  \xi, \eta\in \mathbb{R}^n\backslash \{0\}
\end{align}\[\]
hold provided $H \in C^1(\mathbb{R}^n\backslash \{0\})$ (see \cite{BC}, subsection 2.2). We also notice that \eqref{2.1.3} and \eqref{2.1.4} imply that
\begin{align}\label{2.1.5}
\nabla_{\xi}H(\xi)\cdot\eta  \leq H(\xi)\ \ \forall, \eta\in \mathbb{R}^n\backslash \{0\}.
\end{align}\[\]
Furthermore, the map $H\nabla_{\xi}H$ is invertible with
\vskip0.1cm
\begin{align}\label{2.1.6}
H\nabla_{\xi}H = (H_0\nabla_{\xi}H_0)^{-1}.
\end{align}\[\]
From \eqref{2.1.4} and the homogeneity of $H_0$, \eqref{2.1.6} is equivalent to
\vskip0.1cm
\begin{align}\label{2.1.7}
H(\xi)\nabla_{\eta}H_0(\nabla_{\xi}H(\xi)) = \xi.
\end{align}\[\]
\par Sometimes we write
\vskip0.1cm
\begin{align*}
\Delta^{H}_{p}u = -\mathrm{div}\;(a(\nabla u)),
\end{align*}\[\]
in the form of divergence, where $\Delta^{H}_{p}$ is called the \emph{Finsler $p$-Laplace}(or anisotropic $p$-Laplace) operator and $a(\nabla u)$ is given by \eqref{1.4}. More precisely, \eqref{1.6} reads as
\vskip0.1cm
\begin{align}\label{2.1.8}
-\Delta^{H}_{p}u = u^{p^{*}-1} ,
\end{align}\[\]
where
\begin{align*}
p^{*} = \frac{np}{n-p}.
\end{align*}

\par The following Lemma is a refinement for property of $H$ operator. We omit its proof which is contained in \cite{CFV}.
\begin{lemma}\label{lem 2.1}
Assume that $H$ in $C^2( \mathbb{R}^n\backslash \{0\})$, it holds that
\begin{enumerate}
  \item $\sum\limits_{i=1}^{n}H_{i}(\xi)\xi_{i} = H(\xi)$,\\
  \item $\sum\limits_{i=1}^{n}H_{ij}(\xi)\xi_{i} = 0$,\\
  \item $H_{ij}(t\xi) = \frac{1}{t}H_{ij}(\xi)$.
\end{enumerate}
\end{lemma}
\vskip0.1cm

The regularity theory for Sobolev equation in divergence form, modeled upon the Laplacian, $p$-Laplacian, and anisotropic Laplacian, have extensively been developed in the past years (see \cite{Caffarelli}, \cite{CFR}, \cite{Dibenedetto},  \cite{DM}, \cite{Iwaniec}, \cite{Lewis}, \cite{Lieberman}, \cite{LY}, \cite{LZ}, \cite{Tolksdorf}, \cite{Uhlenbeck} and the references therein). If a more general proof for regularity of anisotropic equation is desired, we recommend that readers refer to Reference \cite{ACF}. We present some results regarding regularity of anisotropic equation here just for completeness and convenience of readers. Notice that Einstein summation convention of summation is used throughout the paper, we will omit the sum sign below.
\medskip
\vskip0.1cm

\noindent\textbf{2.2 Regularity of solutions of critical anisotropic Sobolev equation.}
\medskip
\begin{lemma}[See \cite{ACF}]\label{lem 2.2}
Let $u \in W^{1,p}_{loc}(\Omega)$ be a local weak solution of the equation
\vskip0.1cm
\begin{align}\label{2.2.1}
-\mathrm{div}\;(a(\nabla u)) = f,
\end{align}\[\]
with $f \in L^{q}_{loc}(\Omega)$ and $q$ satisfies
\begin{equation}\label{2.2.2}
q=\begin{cases}
2& p \geq \frac{2n}{n+2}\\
(p^{*})',&1< p < \frac{2n}{n+2}.
\end{cases}\end{equation}\[\]
Then $a(\nabla u)$ belongs to $H^{1}_{loc}(\Omega)$.
\end{lemma}

\begin{lemma}[See \cite{ACF}]\label{lem 2.3}
Let $u\in W^{1,p}_{loc}(\Omega)$ be a local weak solution of the equation
\vskip0.1cm
\begin{align}\label{2.2.1}
-\mathrm{div}\;(a(\nabla u)) = f,
\end{align}\[\]
where $f \in L^{r}_{loc}(\Omega)$ with $r>n$. Then $u \in H^2_{loc}(\Omega)\cap C^{1,\beta}_{loc}(\Omega)$ for $\beta \in (0,1)$ depending only on $n,p,r$ and $H$.
\end{lemma}
\medskip
\vskip0.1cm

\section{A vital integral inequality on vector fields}
In this section, we need some preliminaries before proving Theorem \ref{thm 1.1}. More precisely, the vital integral inequality \eqref{3.3.2} plays a key role in proving Theorem \ref{thm 1.1}. Hence our main goal in this section is to prove the vital integral inequality \eqref{3.3.2}. Before presenting \eqref{3.3.2}, we first define vector fields and show some lemmas that we need.

\medskip
\vskip0.1cm
\noindent\textbf{3.1 Definition of vector fields.}
\medskip
Letting $u > 0$ be any weak solution of \eqref{1.6}, and $u$ satisfies
\begin{equation}\label{3.1.1}
u \in  C^{1,\tau}_{loc}(\mathbb{R}^n)
\end{equation}
in the previous Lemma \ref{lem 2.3}, one could immediately deduce that

\begin{equation}\label{3.1.2}
H^{p-1}(\nabla v)\nabla H(\nabla v) \in W^{1,2}_{loc}(\mathbb{R}^n).
\end{equation}
from Lemma \ref{lem 2.1}.
\medskip

Now we introduce the following vector fields
\begin{equation*}
a^i = H^{p-1}(\nabla u)H_{i}(\nabla u),
\end{equation*}
and
\begin{equation*}
W_{ij} = a^i  ,_{j} - \frac{a^i u_j}{\omega(u)}  - \frac{1}{n}\left( \Delta^{H}_{p}u-\frac{H^{p}(\nabla u)}{\omega(u)}\right)\delta_{ij},
\end{equation*}
where $W_{ij}$ is trace free tensor and $\omega(u)$ could be determined later in Remark \ref{rem 3.2}. With the help of \eqref{3.1.1} and \eqref{3.1.2} , $a^i \in L^{\infty}_{loc}(\mathbb{R}^n)$ and $W_{ij}\in L^2_{loc}(\mathbb{R}^n)$. Denote $p_{*} = \frac{p(n-1)}{n-p}$.
\medskip
\par Recalling the definition of $a^i$, we have

\begin{equation*}
\Delta_{p}^{H}u = a^j_{,j}    \;\;\;\; \text{in} \;\;\mathbb{R}^n,
\end{equation*}\[\]
in the weak sense, that is

\begin{equation}\label{3.1.3}
-\int_{\mathbb{R}^n}a^j\psi_j= \int_{\mathbb{R}^n}\Delta_{p}^{H}u\;\psi.
\end{equation}\[\]

\vskip0.1cm
\noindent\textbf{3.2 Anisotropic differential identities.}
The following Lemma \ref{lem 3.1} can be found in Reference \cite{LWY}, and here we only state the results for brevity. For detailed proofs, we refer the reader to Reference \cite{LWY}. In this subsection, we provide detailed computations only for certain key points.

\begin{lemma}\label{lem 3.1}
With the notations as in above, then we have
\begin{enumerate}
  \item $a^i_{,j}u_i = (p-1)a^k u_{kj}$,\\
  \item $W_{ij}u_i = (p-1)a^k u_{kj} - \frac{n-1}{n}\frac{H^{p}(\nabla u)u_j}{\omega(u)}+ \frac{1}{n}u^{p^{*}-1}u_j$,\\
  \item $W_{ij,i} = \frac{1}{p_{*}-1}\frac{W_{ij}u_i}{\omega(u)} -\frac{u^{p^{*}-1}u_j}{\omega(u)} - \frac{n-p}{n}\frac{a^i u_{ij}}{\omega(u)}+\frac{n-1}{n}\frac{\omega'(u)H^{p}(\nabla u)}{\omega(u)}$.
\end{enumerate}
\end{lemma}

\begin{remark}\label{rem 3.2}
Recalling the definition of $W_{ij}$, we obtain
\begin{align}\label{3.2.1}
\begin{split}
(W_{ij}a^{j})_{,i} &= W_{ij,i}a^{j} + W_{ij}a^{j}_{,i}
\\& = W_{ij,i}a^{j} + W_{ij}\left(W_{ji} + \frac{a^i u_j}{\omega(u)} + \frac{1}{n}( \Delta^{H}_{p}u-\frac{H^{p}(\nabla u)}{\omega(u)})\delta_{ij}\right)
\\& = W_{ij,i}a^{j} + W_{ij}W_{ji} + \frac{W_{ij}a^i u_j}{\omega(u)} + \frac{W_{ij}}{n}\Delta^{H}_{p}u-\frac{W_{ij}}{n}\frac{H^{p}(\nabla u)}{\omega(u)}\delta_{ij}
\\& = W_{ij}W_{ji} + (1- \frac{1}{p_{*}-1}) \frac{W_{ij}a^{j}u_i}{\omega(u)} -\frac{n-1}{n}\left((p^{*}-1)u^{p^{*}-2}- \frac{p^{*}-1}{p_{*}-1}\frac{u^{p^{*}-1}}{\omega(u)}\right)H^{p}(\nabla u)\\&+\frac{n-1}{n}(\omega'(u) - \frac{1}{p_{*}-1})\frac{H^{2p}(\nabla u)}{\omega^{2}(u)}.
\end{split}
\end{align}
If we take $\omega(u) = \frac{u}{p_{*}-1}$, then we get $(W_{ij}a^{j})_{,i} = W_{ij}W_{ji}$. Hence we will replace $\omega(u)$ by $\frac{u}{p_{*}-1}$ below for calculate.
\end{remark}

\vskip0.1cm
\begin{remark}\label{rem 3.3}
Using the same way as \eqref{3.2.1}, we define the function $g = u^{\alpha -1}H^{p}(\nabla u) + \beta u^{\alpha}\Delta_{p}^{H}u$ and obtain that
\begin{align}\label{3.2.2}
\begin{split}
g_{i} &= (\alpha -1)u^{\alpha -2}H^{p}(\nabla u) +p u^{\alpha -1}H^{p-1}(\nabla u)H_{k}(\nabla u)u_{ki} + \alpha\beta u^{\alpha-1}u_i\Delta_{p}^{H}u + \beta u^{\alpha}(\Delta_{p}^{H}u)_{i}
\\& = \left(\alpha -1 + \frac{p(n-1)(p_{*}-1)}{n(p-1)}\right)u^{\alpha -2}H^{p}(\nabla u)u_i + \frac{p}{p-1}u^{\alpha-1}W_{ij}u_i \\&+ (\frac{p}{n(p-1)}+ \alpha\beta + \beta(p^{*}-1))u^{\alpha-1}u_i\Delta^{H}_{p}u.
\end{split}
\end{align}
If we take $\alpha = -\frac{n(p-1)}{n-p}$ and $\beta = -\frac{n-p}{n(p-1)}$, then we get $g_{i} = \frac{p}{p-1}u^{-p_{*}}W_{ij}u_i$ and $g = u^{-p_{*}}H^{p}(\nabla u) + \frac{n-p}{n(p-1)} u^{\frac{p}{n-p}}$.
\end{remark}

\begin{lemma}\label{lem 3.4}
\begin{align}
\begin{split}
[\tau(u)]^{-1}(\tau(u)W_{ij}a^j)_{,i} &= W_{ij}W_{ji} + (\frac{\tau'(u)}{\tau(u)} + \frac{p_{*}-2}{p_{*}-1}\frac{1}{\omega(u)})W_{ij}a^ju_{i} \\&- \frac{n-1}{n}((p^{*}-1)u^{p^{*}-2}- \frac{p^{*}-1}{p_{*}-1}\frac{u^{p^{*}-1}}{\omega(u)}  )H^{p}(\nabla u) \\
&+ \frac{n-1}{n}(\omega'(u) - \frac{1}{p_{*}-1})\frac{H^{2p}(\nabla u)}{\omega^{2}(u)}.
\end{split}
\end{align}\[\]
\end{lemma}
\begin{proof}
Combining with Remark \ref{rem 3.2} and the definition of $a^i$, and using the statement (3) of Lemma \ref{lem 2.1}, we obtain
\begin{align}\label{3.2.3}
\begin{split}
[\tau(u)]^{-1}(\tau(u)W_{ij}a^j)_{,i}
&  = \frac{\tau'(u)u_i}{\tau(u)} W_{ij}a^j + (W_{ij}a^j)_{,i}                                                  \\
&  = \frac{\tau'(u)u_i}{\tau(u)} W_{ij}a^j +     (W_{ij})_{,i} a^j +   W_{ij}(a^j)_{,i}\\
&  =  W_{ij}W_{ji} + (\frac{\tau'(u)}{\tau(u)} + \frac{p_{*}-2}{p_{*}-1}\frac{1}{\omega(u)})W_{ij}\overrightarrow{\mathbf{v}}^ju_{i} \\&- \frac{n-1}{n}((p^{*}-1)u^{p^{*}-2}- \frac{p^{*}-1}{p_{*}-1}\frac{u^{p^{*}-1}}{\omega(u)})H^{p}(\nabla u) \\
&+ \frac{n-1}{n}(\omega'(u) - \frac{1}{p_{*}-1})\frac{H^{2p}(\nabla u)}{\omega^{2}(u)} . \\
\end{split}
\end{align}
\end{proof}

\vskip0.1cm
\begin{remark}\label{rem 3.5}
If we choose $\tau(u)  = u^{2-p_{*}}$ and $\omega(u) = \frac{p_{*}-1}{u}$, then we obtain that
\begin{align}\label{3.2.4}
\begin{split}
(u^{2-p_{*}}W_{ij}a^j)_{,i} = u^{2-p_{*}}W_{ij}W_{ji}.
\end{split}
\end{align}
\end{remark}

\medskip
\vskip0.1cm
\noindent\textbf{3.3 The vital differential inequality.}
\medskip
The regularity of identities is an important topic in the study of partial differential equations. However, in this paper, we do not focus on proving the regularity of such identities. There is already a vast amount of literature on the regularity of various non-homogeneous equations such as  \cite{CFR}, \cite{CFV},  \cite{CG}, \cite{LY}, \cite{LZ}, \cite{Serrin}, \cite{SZ} \cite{Zhou} and its refenences. Therefore, in this subsection, we aim to prove the vital differential inequality \eqref{3.3.10} which plays an important role in proving Theorem \ref{thm 1.1} and assume the relevant regularity results hold. For more details, we refer readers to Zhou's work in Reference \cite{Zhou}, where the process and results are presented. Let $\rho$ be a smooth cut-off function satisfying:
\begin{equation}\label{rho}
\\\begin{cases}
\rho \equiv 1 &  in \;\; B_R, \\
0 \leq \rho \leq 1 & in \;\; B_{2R},\\
\rho \equiv 0 & in \;\;\mathbb{R}^n \backslash B_{2R},\\
|\nabla \rho| \lesssim \frac{1}{R} & in \;\;\mathbb{R}^n,
\\\end{cases}
\end{equation}
where and in the sequel. Moreover we use $"\lesssim"$, $"\backsimeq"$ to replace $"\leq"$, $"="$, etc., to drop out some positive constants independent of $R$ and $v$.
\vskip0.1cm
\begin{lemma}\label{lem 3.6}
Let $n \geq 2$, $1 < p <n$ and $ u \in W^{2,2}_{loc}(\mathbb{R}^n) \cap  C^{1,\alpha}_{loc}(\mathbb{R}^n)$ be a positive, weak solution of \eqref{1.6}, $W$ be the $n\times n$ square matrix whose elements are denoted by $\{W_{ij}\}$. Then for all $i,j,k \in \mathcal{I}$ we have
\begin{align*}
\begin{split}
0 \leq \sum\limits_{i,j,k=1}^{n}W_{ij}a^{j}W_{ki}u_k \leq \Lambda H^{p}(\nabla u)\sum\limits_{i,j=1}^{n}W_{ij}W_{ji},
\end{split}
\end{align*}
where $\Lambda = \lambda\max\limits_{\substack{\xi \in \mathbb{R}^n\\1 \leq i, j \leq n}}\left\{\frac{(\nabla^{2}_{ij}H^{p}(\xi))|\xi|^{2}} {p(p-1)H^{p}(\xi)}\right\}$.
\end{lemma}
\begin{proof}
We first observe that $(a^{i}(\nabla v))_{,j} = \{(H^{p-1}(\nabla v)H_{i}(\nabla v))_{j}\}_{n\times n} = AC$, with $C$ is Hessian matrix of $u$ and $A = (p-1)H^{p-2}(\nabla v)\nabla H(\nabla v)\bigotimes \nabla H(\nabla v) + H^{p-1}(\nabla v)\nabla^2H(\nabla v)$. Since $H^2$ is uniformly convex, Hessian matrix of $H^2$ is positive definite and we obtain that the matrix $A$ is positive definite and symmetric. Then we can rewrite $W= AB - \frac{1}{n}T_r(AB)I_n$ and $B= C - \frac{\nabla v\bigotimes \nabla v}{(p-1)\omega(v)}$. Obviously, $W =AB $ if $i\neq j$. $H(\nabla v)$ written as $H$ and $a_j(\nabla v)$ written as $a_j$ for convenience and careful computation gives
\begin{align}\label{3.3.1}
\begin{split}
\sum\limits_{i,j,k=1}^{n}W_{ij}a^{j}W_{ki}u_k &= \sum\limits_{i=j=k}W_{ii}a^{i}W_{ii}u_i +  \sum\limits_{i\neq j\neq k}W_{ij}a^{j}W_{ki}u_k + \sum\limits_{i=j\neq k}W_{ii}a^{i}W_{ki}u_k \\& + \sum\limits_{i=k\neq j}W_{ij}a^{j}W_{ii}u_i + \sum\limits_{i\neq j = k}W_{ij}a^{j}W_{ji}u_j                                  \\
 &= \sum\limits_{i=j=k}W_{ii}a^{i}W_{ii}u_i +  \sum\limits_{i\neq j\neq k}A_{im}B_{mj}a_{j}A_{kt}B_{ti}u_k + \sum\limits_{i=j\neq k}W_{ii}a_{i}A_{km}B_{mi}u_k
 \\&+ \sum\limits_{i=k\neq j}A_{im}B_{mj}a^{j}W_{ii}u_i + \sum\limits_{i\neq j = k}A_{im}B_{mj}a^{j}A_{jt}B_{ti}u_j\\
 &= \sum\limits_{i=j=k}W_{ii}^{2}a^{i}u_i +  \sum\limits_{i\neq j\neq k}A_{ii}B_{ij}a^{j}A_{kk}B_{ki}u_k + \sum\limits_{i=j\neq k}W_{ii}a^{i}A_{kk}B_{ki}u_k
 \\&+ \sum\limits_{i=k\neq j}A_{ii}B_{ij}a^{j}W_{ii}u_i + \sum\limits_{i\neq j = k}A_{ii}B_{ij}a^{j}A_{jj}B_{ji}u_j\\
 &= \frac{1}{p-1}\sum\limits_{i=j=k}W_{ii}^{2}A_{ii}u_{i}u_{i} + \frac{1}{p-1}\sum\limits_{i\neq j\neq k}A_{ii}B_{ij}A_{jj}u_{j}A_{kk}B_{ki}u_k \\&+ \frac{1}{p-1}\sum\limits_{i=j\neq k}W_{ii}A_{ii}u_{i}A_{kk}B_{ki}u_k + \frac{1}{p-1}\sum\limits_{i=k\neq j}A_{ii}B_{ij}A_{jj}u_{j}W_{ii}u_i  \\&+ \frac{1}{p-1}\sum\limits_{i\neq j = k}A_{ii}B_{ij}A_{jj}u_{j}A_{jj}B_{ji}u_j\\
 &= I_{1} +I_{2} + I_{3} + I_{4} + I_{5}.
\end{split}
\end{align}
For the second term $I_2$, we observe that
\begin{align}\label{3.3.2}
\begin{split}
 \frac{1}{p-1}\sum\limits_{i\neq j\neq k}A_{ii}B_{ij}A_{jj}u_{j}A_{kk}B_{ki}u_k &\leq \frac{1}{2(p-1)}\sum\limits_{i\neq j\neq k}\left(A_{ii}B_{ij}A_{jj}^2u_{k}^2 + A_{ii}A_{jj}^2B_{ki}^2u_{j}^2\right)
 \\& = \frac{1}{p-1}\sum\limits_{i\neq j\neq k}A_{ii}B_{ij}A_{jj}^2u_{k}^2
 \\& = \frac{1}{p-1}\sum\limits_{i\neq j}A_{ii}B_{ij}A_{jj}^2(|\nabla u|^2 - u_{i}^2 - u_{j}^2)
\end{split}
\end{align}
 For the third and fourth term $I_3 , I_4$, we observe that
\begin{align}\label{3.3.3}
\begin{split}
 &\frac{1}{p-1}\sum\limits_{i=j\neq k}W_{ii}A_{ii}u_{i}A_{kk}B_{ki}u_k + \frac{1}{p-1}\sum\limits_{i=k\neq j}A_{ii}B_{ij}A_{jj}u_{j}W_{ii}u_i \\&= \frac{2}{p-1}\sum\limits_{i\neq j}A_{ii}B_{ij}A_{jj}u_{j}W_{ii}u_i
 \\& \leq \frac{1}{p-1}\sum\limits_{i\neq j} \left(A_{ii}A_{jj}^2B_{ij}^2u_{i}^2 + A_{ii}W_{ii}^2u_{j}^2 \right)
 \\& =  \frac{1}{p-1} \sum\limits_{i\neq j} A_{ii}A_{jj}^2B_{ij}^2u_{i}^2 +  \frac{1}{p-1}\sum\limits_{i}A_{ii}W_{ii}^2 (|\nabla u|^2 - u_{i}^2)
 \\& = I_6 + I_7.
\end{split}
\end{align}
 For the fifth term $I_5$
\begin{align}\label{3.3.4}
\begin{split}
\frac{1}{p-1}\sum\limits_{i\neq j = k}A_{ii}B_{ij}A_{jj}u_{j}A_{jj}B_{ji}u_j  = \frac{1}{p-1}\sum\limits_{i\neq j}A_{ii}A_{jj}^2B_{ij}^2u_{j}^2.
\end{split}
\end{align}
Furthermore, we compute the above inequalities and yields that
\begin{align}\label{3.3.5}
\begin{split}
I_1 + I_7 =  \frac{1}{p-1}\sum\limits_{i=1}^{n}A_{ii}W_{ii}^2|\nabla v|^2,
\end{split}
\end{align}
and
\begin{align}\label{3.3.6}
\begin{split}
I_2 + I_3 + I_4 +I_5 +I_6 &\leq  \frac{1}{p-1} \sum\limits_{i\neq j} A_{ii}A_{jj}^2B_{ij}^2|\nabla u|^2
\\& =  \frac{1}{p-1} \sum\limits_{i < j} A_{ii}A_{jj}(A_{ii} +A_{jj})B_{ij}^2|\nabla u|^2.
\end{split}
\end{align}
Hence, combining with above inequalities, we obtain that

\begin{align}\label{3.3.7}
\begin{split}
\sum\limits_{i,j,k=1}^{n}W_{ij}a^{j}W_{ki}u_k &\leq \frac{1}{p-1}\sum\limits_{i=1}^{n}A_{ii}W_{ii}^2|\nabla u|^2 +  \frac{1}{p-1} \sum\limits_{i < j} A_{ii}A_{jj}(A_{ii} +A_{jj})B_{ij}^2|\nabla u|^2
\\& \leq \Lambda H^{p}\sum\limits_{i=1}^{n}W_{ii}^2 + \Lambda H^{p} \sum\limits_{i < j}A_{ii}B_{ij}A_{jj}B_{ji}
\\&  = \Lambda \left(H^{p} \sum\limits_{i= j}W_{ii}^2  + H^{p} \sum\limits_{i\neq j}W_{ij}W_{ji}\right)
\\& = \Lambda H^{p}\sum\limits_{i,j=1}^{n}W_{ij}W_{ji}.
\end{split}
\end{align}

\par Finally, we will prove that $ \sum\limits_{i,j,k=1}^{n}W_{ij}a^{j}W_{ki}u_k $ is non-negative. We define matrix $K = a^ju_k$, then one should be noted that matrix $K$ is the idempotent matrix. Since eigenvalues of idempotent matrices $K$ are 0 or 1 and the rank of $K$ is 1, there exists an invertible matrix $T$ such that $T^{-1}KT$ is diagonal matrix with eigenvalues $\lambda_j$ is equal to 1 for fixed $j$. We may assume that $j=1$ and careful computation gives

\begin{align}\label{3.3.8}
\begin{split}
\sum\limits_{i,j,k=1}^{n}W_{ij}K_{jt}W_{ti} &= \sum\limits_{i}^{n}W_{ij}K_{jj}W_{ji}
\\&= \sum\limits_{i=1}^{n}W_{i1}W_{1i}
\\& = W_{11}^2 + \sum\limits_{i=2}^{n}W_{i1}W_{1i}
\end{split}
\end{align}
For the second term, there exists an orthogonal matrix $T$ such that $T^{-1}AT$ is a diagonal matrix. Define $\widetilde{A} = T^{-1}AT$ and $\widetilde{B} = T^{-1}BT$, where $\widetilde{A} = \{\widetilde{a_{ij}}\}$ and $\widetilde{B} = \{\widetilde{b_{ij}}\}$  are diagonal matrices and careful computation gives we obtain that
\begin{align}\label{3.3.9}
\begin{split}
 T_r\{W^2\} &=  T_r\{ABAB\}\\
 &= T_r\{T^{-1}ABABT\}\\
 &= T_r\{T^{-1}ATT^{-1}BTT^{-1}ATT^{-1}BT\}\\
 &= T_r\{\widetilde{A}\widetilde{B}\widetilde{A}\widetilde{B}\}\\
 &= \sum\limits_{i,j,k,l=1}^{n}\widetilde{a_{ij}}\widetilde{b_{jk}}\widetilde{a_{kl}}\widetilde{b_{li}}\\
 &= \sum\limits_{i,k=1}^{n}\widetilde{a_{ii}}\widetilde{a_{kk}}\widetilde{b_{ik}}^{2}\geq 0,\\
\end{split}
\end{align}
since $A$ is positive definite. Hence we finish the proof.
\end{proof}

\medskip
\vskip0.1cm

\begin{remark}
We remark that $\Lambda \geq 1$ holds for the reason that
\begin{align*}
\begin{split}
p(p-1)H^{p}(\nabla u) = \sum\limits^{n}_{i,j=1}\nabla^{2}_{ij}H^{p}(\xi)\xi_i\xi_j \leq |\xi|^2 \lambda\max\nabla^{2}_{ij}H^{p}(\xi).
\end{split}
\end{align*}
\end{remark}

\medskip
\vskip0.1cm

\begin{remark}
When  $H(\xi) = |\xi|$, then $\Lambda = 1$. Our theorem closely aligns with the result obtained by V\'{e}tois in \cite{Vet24}. Nevertheless, under this condition, a rotational transformation combined with the trace-free property of W leads to $\Lambda = \frac{n-1}{n}$, in which case this results coincide with those reported by Sun-Wang in \cite{SW} and the range of $p$ will become $p_n < p <n$ where
\begin{equation*}
p_n=
\begin{cases}
\frac{n^2}{3n-2}&\mathrm{if}\;\;\;n=2,3,4,\\
\frac{n^2+2}{3n}&\mathrm{if}\;\;\;n\geq 5 .
\end{cases}
\end{equation*}

\end{remark}

\begin{lemma}\label{lem 3.9}
For $0 < m < \frac{p-1}{p\Lambda}$ and $\varepsilon > 0$, we have
\begin{align}\label{3.3.10}
\begin{split}
(g^{-m}u^{2-p_{*}}W_{ij}a^{j})_{,i} \geq  \varepsilon g^{-m}u^{2-p_{*}}W_{ij}W_{ji},
\end{split}
\end{align}\[\]
where $ \varepsilon = 1-\frac{pm\Lambda}{p-1}$.
\end{lemma}
\begin{proof}
\begin{align*}
\begin{split}
(g^{-m}u^{2-p_{*}}W_{ij}a^{j})_{,i} &= -mg^{-m-1}g_{i}u^{2-p_{*}}W_{ij}a^{j} + g^{-m}u^{2-p_{*}}W_{ij}W_{ji}\\
&= -\frac{pm}{p-1}g^{-m-1}u^{2-2p_{*}}W_{ij}a^{j}W_{ki}u_i + g^{-m}u^{2-p_{*}}W_{ij}W_{ji}\\
& \geq  \left(1-\frac{pm\Lambda}{p-1}\right)g^{-m}u^{2-p_{*}}W_{ij}W_{ji}\\
& = \varepsilon g^{-m}u^{2-p_{*}}W_{ij}W_{ji}.
\end{split}
\end{align*}\[\]
\end{proof}

\begin{lemma}\label{lem 3.8}
If $\alpha < 0$, then we have
\begin{align}\label{3.3.11}
\begin{split}
&\int_{\mathbb{R}^n}u^{\alpha-1}g^{\beta}H^{p}(\nabla u)\rho^{\gamma} + \int_{\mathbb{R}^n}u^{\alpha-1+p^{*}}g^{\beta}\rho^{\gamma} \\&\lesssim \int_{\mathbb{R}^n}u^{\alpha}g^{\beta}H^{p-1}(\nabla u)\rho^{\gamma-1}|\nabla \rho| + \int_{\mathbb{R}^n}u^{\alpha-p_{*}}g^{\beta-1}H^{p}(\nabla u)|W_{ij}|\rho^{\gamma}.
\end{split}
\end{align}\[\]
\end{lemma}
\begin{proof}
Through a straightforward calculation, we can readily obtain that
\begin{align*}
\begin{split}
(u^{\alpha}g^{\beta}a^{i})_{,i} &= \alpha u^{\alpha-1}u_{i}g^{\beta}a^{i} + \beta u^{\alpha}g^{\beta-1}g_{i}a^{i} + u^{\alpha}g^{\beta}a^{i}_{,i}\\
& = \alpha u^{\alpha-1}g^{\beta}H^{p-1}(\nabla u) + \beta u^{\alpha}g^{\beta-1}(\frac{p}{p-1}u^{-p_{*}}W_{ij}u_i)a^{i} - u^{\alpha}g^{\beta}u^{p^{*}-1}.
\end{split}
\end{align*}\[\]
Testing both side of equation by the test function $\rho^{\gamma}$, integrating and applying integration by parts, we obtain
\begin{align*}
\begin{split}
&\int_{\mathbb{R}^n}u^{\alpha-1}g^{\beta}H^{p}(\nabla u)\rho^{\gamma} + \int_{\mathbb{R}^n}u^{\alpha-1+p^{*}}g^{\beta}\rho^{\gamma} \\&\lesssim \int_{\mathbb{R}^n}u^{\alpha}g^{\beta}H^{p-1}(\nabla u)\rho^{\gamma-1}|\nabla \rho| + \int_{\mathbb{R}^n}u^{\alpha-p_{*}}g^{\beta-1}H^{p}(\nabla u)|W_{ij}|\rho^{\gamma}.
\end{split}
\end{align*}
\end{proof}

\medskip
\vskip0.1cm
\noindent\textbf{3.4 Asymptotic estimates on bounded region.} The main goal of this subsection is to prove asymptotic estimates below. We first prove Lemma \ref{lem 3.11}. Corollary \ref{cor 3.12} and Corollary \ref{cor 3.13} are two important generalizations of Lemma \ref{lem 3.11}.

\begin{lemma}\label{lem 3.11}
For $-p \leq q < -1$, we have
\begin{equation}\begin{split}\label{3.4.1}
\int_{B_R}u^{q}H^{p}(\nabla u)  \lesssim  R^{n-\frac{p(p^{*}+q)}{p^{*}-p}},
\end{split}\end{equation}
and
\begin{equation}\begin{split}\label{3.4.2}
\int_{B_R}u^{q}  \lesssim  R^{n-\frac{p(p^{*}+q)}{p^{*}-p}}.
\end{split}\end{equation}
\end{lemma}

\begin{proof}
Since $u$ is the solution of \eqref{1.6} in weak sense, we have
\begin{equation}\label{3.4.3}
-\int_{\mathbb{R}^n}a^{j}\psi_j= \int_{\mathbb{R}^n} \Delta^{H}_{p}u \;\psi
\end{equation}

Replacing $\psi$ by $u^{1+q}\psi $ in \eqref{3.4.3}, then we consider the term on the left side of \eqref{3.4.3} to derive
\begin{align}\label{3.4.4}
\begin{split}
 -\int_{\mathbb{R}^n}a^{j}(u^{1+q}\psi)_j &= -\int_{\mathbb{R}^n} (1+q)a^{j}u^{q}u_j\psi -\int_{\mathbb{R}^n} a^{j}u^{1+q}\psi_j
 \\& = -\int_{\mathbb{R}^n} (1+q) H^{p}(\nabla u)u^{q}\psi -\int_{\mathbb{R}^n} a^ju^{1+q}\psi_j.
\end{split}
\end{align}
From the term on the right side of \eqref{3.4.3}, we obtain
\begin{align}\label{3.4.5}
\begin{split}
\int_{\mathbb{R}^n} \Delta^{H}_{p}u\;u^{1+q}\psi &= -\int_{\mathbb{R}^n} u^{p^{*}-1}u^{1+q}\psi
 \\ &= -\int_{\mathbb{R}^n}  u^{p^{*}+q}\psi.
\end{split}
\end{align}
Combining with \eqref{3.4.4} and \eqref{3.4.5}, we get
\begin{equation}\label{3.4.6}
-\int_{\mathbb{R}^n} (1+q) H^{p}(\nabla u)u^{q}\psi   + \int_{\mathbb{R}^n}  u^{p^{*}+q}\psi = \int_{\mathbb{R}^n} a^j u^{1+q}\psi_j.
\end{equation}
Next, let $\theta > 0$ be a constant big enough and $\rho$ be the cut-off function as \eqref{rho}. Using \eqref{3.4.6} with $\psi = \rho^{\theta}$ we have
\begin{align}\label{3.4.7}
\begin{split}
-\int_{\mathbb{R}^n} (1+q) H^{p}(\nabla u)u^{q}\rho^{\theta}   + \int_{\mathbb{R}^n}  u^{p^{*}+q}\rho^{\theta}  &= \theta\int_{\mathbb{R}^n} a^j u^{1+q}\rho^{\theta-1}\rho_j.
\end{split}
\end{align}
Since $H(\nabla v)\in C^{1,\tau}_{loc}(\mathbb{R}^n)$, there exists a constant $M_1 > 0$ such that $H_{j}(\nabla v) \leq M_1$. Moreover, since $|a^j\rho_j| \lesssim \frac{M_1}{R}H^{p-1}(\nabla v)$, we derive that
\begin{align}\label{3.4.8}
\begin{split}
\theta\int_{\mathbb{R}^n}a^j u^{1+q}\rho^{\theta-1}\rho_j &\lesssim \frac{M_1}{R}\int_{R^n} u^{1+q}H^{p-1}(\nabla u)\rho^{\theta-1}
\\&\leq  \varepsilon\int_{R^n}\rho^{\theta}u^{q} H^{p}(\nabla u)  + \frac{M_{1}^p}{\varepsilon^{p-1}R^p}\int_{R^n}u^{p+q}\rho^{\theta-p},
\end{split}
\end{align}
using the $\varepsilon$-Young's inequality with exponent pair $(\frac{p}{p-1},p)$. Then it follows that
\begin{align}\label{3.4.9}
  -\int_{\mathbb{R}^n} (1+q) H^{p}(\nabla u)u^{q}\rho^{\theta}   + \int_{\mathbb{R}^n}  u^{p^{*}+q}\rho^{\theta}  \leq \varepsilon\int_{\mathbb{R}^n} \rho^{\theta}u^{q} H^{p}(\nabla u)  + \frac{M_{1}^p}{\varepsilon^{p-1}R^p}\int_{R^n}u^{p+q}\rho^{\theta-p}.
\end{align}
Using $\varepsilon$-Young's inequality with exponent $(\frac{p^{*}+q}{p+q},\frac{p^{*}+q}{p^{*}-p})$, then we derive that
\begin{align}\label{3.4.10}
\begin{split}
   \frac{M_{1}^p}{\varepsilon^{p-1}R^p}\int_{R^n}u^{p+q}\rho^{\theta-p} \leq \varepsilon \int_{R^n}u^{p^{*}+q}\rho^{\theta} + \frac{(M_1)^{\frac{p(p^{*}+q)}{p^{*}-p}}}{\varepsilon^{\frac{pp^{*}+pq-p^{*}+p}{p^{*}-p}}} R^{\frac{-p(p^{*}+q)}{p^{*}-p}} \int_{R^n} \rho^{\frac{\theta p^{*}-pp^{*}-pq-\theta p}{p^{*}-p}}.
\end{split}
\end{align}
Inserting \eqref{3.4.10} into \eqref{3.4.9} yields
\begin{align*}
\begin{split}
&-\int_{\mathbb{R}^n} (1+q) H^{p}(\nabla u)u^{q}\rho^{\theta}   + \int_{\mathbb{R}^n}  u^{p^{*}+q}\rho^{\theta}  \\&\lesssim \varepsilon\int_{\mathbb{R}^n}\rho^{\theta}u^{q} H^{p}(\nabla u)  + \varepsilon \int_{R^n}u^{p^{*}+q}\rho^{\theta} + \frac{(M_1)^{\frac{p(p^{*}+q)}{p^{*}-p}}}{\varepsilon^{\frac{pp^{*}+pq-p^{*}+p}{p^{*}-p}}} R^{\frac{-p(p^{*}+q)}{p^{*}-p}} \int_{R^n} \rho^{\frac{\theta p^{*}-pp^{*}-pq-\theta p}{p^{*}-p}}.
\end{split}
\end{align*}
If $-(1+q) > 0$, i.e. $q < -1$, recalling the definition of the test function $\rho$ in \eqref{rho} and taking $\varepsilon > 0$ small enough with $\theta > \frac{p (p^{*}+q)}{p^{*}-p}$, we see that
\begin{align}\label{3.4.11}
\int_{B_R}u^{q}H^{p}(\nabla u) + \int_{B_R}u^{q}  &\lesssim  R^{n-\frac{p(p^{*}+q)}{p^{*}-p}}
\end{align}
This implies \eqref{3.4.1} and \eqref{3.4.2} for $-p < q < -1$. For $q = -p$, a straightforward computation such that \eqref{3.4.1} and \eqref{3.4.2} still valid from \eqref{3.4.9}. Hence we  complete the proof of Lemma \ref{lem 3.9}.
\end{proof}

\begin{corollary}\label{cor 3.12}
Let $n \geq 2$, $1 < p <n$, $0 \leq r \leq p$, $q \leq \frac{n(p-1-r)+p}{n-p}$, $u$ be any weak solution of \eqref{1.6}. Then
\begin{equation*}
\begin{split}
\int_{B_R}u^{q}H^{r}(\nabla u)  \lesssim  R^{n-\frac{n-p}{p}q-\frac{n}{p}r},
\end{split}
\end{equation*}
for $ -r \leq q < \frac{n(p-1-r)+p}{n-p}$, and
\begin{equation*}
\begin{split}
\int_{B_R}u^{q}H^{r}(\nabla u)  \lesssim  R^{n-\frac{n-p}{p-1}q-\frac{n-1}{p-1}r},
\end{split}
\end{equation*}
for $q < -r$.
\end{corollary}

\begin{proof}
Using the conclusion of Lemma \ref{lem 3.8} and Young inequality, for $-r \leq q < \frac{n(p-1-r)+p}{n-p}$, using the Holder inequality we derive that
\begin{align}\label{3.4.18}
\begin{split}
 \int_{\mathbb{R}^n}u^{q}H^{r}(\nabla u) &\leq  \left(\int_{\mathbb{R}^n}u^{q+\sigma(p-r)}H^{r}(\nabla u)\right)^{\frac{r}{p}}\left(\int_{\mathbb{R}^n}u^{q-\sigma r}\right)^{\frac{p-r}{p}}\\
 & \lesssim \left(R^{n-\frac{p(p^{*}+q+\sigma(p-r))}{p^{*}-p}}\right)^{\frac{r}{p}}\cdot \left(R^{n-\frac{p(q-\sigma r)}{p^{*}-p}}\right)^{\frac{p-r}{p}}\\
  &=  R^{n-\frac{n-p}{p}q-\frac{n}{p}r},
\end{split}
\end{align}
where $ \sigma = \max\{\frac{p-q}{p-r},0\}$. Finally, we consider the case where $q < -r$ and combine with \eqref{3.4.18}. Then it follows that
\begin{align}\label{3.4.19}
\begin{split}
 \int_{\mathbb{R}^n}u^{q}H^{r}(\nabla u) &\lesssim R^{-\frac{(n-p)(q+r)}{p-1}}\int_{\mathbb{R}^n}u^{-r}H^{r}(\nabla u)\\
 &\lesssim  R^{n-\frac{n-p}{p-1}q-\frac{n-1}{p-1}r}.
\end{split}
\end{align}
\end{proof}

\medskip
\vskip0.1cm

\begin{corollary}\label{cor 3.13}
Let $n \geq 2$, $1 < p <n$, $\tau \leq 1$ and $\frac{p}{n-p}\tau + \mu < \frac{np-n+p}{n-p}$, $u$ be any weak solution of \eqref{1.6}. Then
\begin{equation}\begin{split}\label{3.4.20}
\int_{B_R}g^{\tau}u^{\mu} \lesssim R^{\max\{n-\tau-\frac{n-p}{p}\mu, n-\frac{n-p}{p-1}\mu\}},
\end{split}\end{equation}
for $ 0 \leq \tau \leq 1$, and
\begin{equation}\begin{split}\label{3.4.21}
\int_{B_R}g^{\tau}u^{\mu}  \lesssim  R^{\max\{n-\tau-\frac{n-p}{p}\mu, n- \frac{p}{p-1}\tau -\frac{n-p}{p-1}\mu\}},
\end{split}\end{equation}
for $\tau < 0$.
\end{corollary}

\begin{proof}
For $ 0 \leq \tau \leq 1$, using $g^{\tau} \leq C_1(u^{-p_{*}\tau}H^{p\tau}(\nabla u) + u^{\frac{n-p}{p}\tau})$ and Corollary \ref{cor 3.12} we obtain that
\begin{align}\label{3.4.22}
\begin{split}
  \int_{B_R}g^{\tau}u^{\mu} &\lesssim  \int_{B_R}(u^{-p_{*}\tau}H^{p\tau}(\nabla u) + u^{\frac{n-p}{p}\tau}) u^{\mu}
   \\&=   \int_{B_R}u^{-p_{*}\tau+\mu}H^{p\tau}(\nabla u) + u^{\tau+\mu}
   \\&\lesssim R^{\max\{n-\tau-\frac{n-p}{p}\mu, n-\frac{n-p}{p-1}\mu\}}.
\end{split}
\end{align}
For $ \tau < 0$, using $g^{\tau} \leq C_2 u^{\frac{p}{n-p}\tau}$ and Corollary \ref{cor 3.12} we obtain that
\begin{align}\label{3.4.23}
\begin{split}
  \int_{B_R}g^{\tau}u^{\mu} &\lesssim  \int_{B_R}u^{\frac{p}{n-p}\tau}\cdot u^{\mu}
   \\&=   \int_{B_R}v^{\frac{p}{n-p}\tau+\mu}
   \\&\lesssim R^{\max\{n-\tau-\frac{n-p}{p}\mu, n- \frac{p}{p-1}\tau- \frac{n-p}{p-1}\mu\}}.
\end{split}
\end{align}
\end{proof}

\section{ Proof of Theorem \ref{thm 1.1}}
In this section, our main effort is to prove the Theorem 1.1 by constructing the correlation between the matrix $W= \{W_{ij}\}$ and vector fields and applying the vital integral inequality \eqref{3.3.11} what we got in Lemma \ref{lem 3.8}. More precisely, we will prove $Tr\{W^{2}\} = 0$ and obtain $W=0$ by the definition of $W$ and careful calculation since $W$ may not be symmetric.

\noindent\textbf{4.2 Proof of Theorem \ref{thm 1.1}}

\begin{proof}
For $p > \frac{1+n\Lambda}{1+2\Lambda}$ and $\Lambda \leq \frac{n-1}{n-2}$. Multiply both sides of \ref{3.3.10} by the test function $\rho^{\gamma}$ and integrate. Applying integration by parts, we obtain that
\begin{align*}
\begin{split}
\int_{\mathbb{R}^n} g^{-m}u^{2-p_{*}}W_{ij}W_{ji}\rho^{\gamma} &\leq C\int_{\mathbb{R}^n}(g^{-m}u^{2-p_{*}}W_{ij}a^{j})_{,i}\rho^{\gamma}\\
& \lesssim \varepsilon_0 \int_{\mathbb{R}^n}g^{-m}u^{2-p_{*}}W_{ij}W_{ji}\rho^{\gamma} +\frac{M_1}{\varepsilon_0} \int_{\mathbb{R}^n}g^{-m}u^{2-p_{*}}H^{2p-2}(\nabla u)\rho^{\gamma-2}|\nabla \rho|^{2}.
\end{split}
\end{align*}
Then we only need to prove that
\begin{align*}
\begin{split}
\int_{\mathbb{R}^n} g^{-m}u^{2-p_{*}}W_{ij}W_{ji}\rho^{\gamma} \lesssim\int_{\mathbb{R}^n}g^{-m}u^{2-p_{*}}H^{2p-2}(\nabla u)\rho^{\gamma-2}|\nabla\rho|^{2}.
\end{split}
\end{align*}\[\]
Taking $m = \frac{p-1}{p\Lambda}-\varepsilon_1$ and combining with $g^{\tau}\leq C\left(u^{p_{*}\tau}H^{p\tau}(\nabla u)+ u^{\frac{p}{n-p}\tau}\right)$, then we have
\begin{align*}
\begin{split}
\int_{\mathbb{R}^n}g^{-m}u^{2-p_{*}}H^{2p-2}(\nabla u)\rho^{\gamma-2}|\nabla \rho|^{2}
&=\int_{\mathbb{R}^n}g^{-\frac{p-1}{p\Lambda}+\varepsilon_1}u^{\frac{2n-p-np}{n-p}}H^{2p-2}(\nabla u)\rho^{\gamma-2}|\nabla \rho|^{2}\\
&\lesssim \int_{\mathbb{R}^n} \left(u^{-p_{*}}H^{p}(\nabla u)\right)^{-\frac{p-1}{p\Lambda}+\varepsilon_1} u^{\frac{2n-p-np}{n-p}}H^{2p-2}(\nabla u)\rho^{\gamma-2}|\nabla \rho|^{2}\\
&\leq  R^{-2}\int_{\mathbb{R}^n}u^{\frac{2n-p-pn}{n-p}+\frac{(n-1)(p-1)}{(n-p)\Lambda}-\frac{p(n-1)}{n-p}\varepsilon_1}H^{\frac{(2p-2)\Lambda-(p-1)}{\Lambda}+p\varepsilon_1}(\nabla u)\rho^{\gamma-2}.
\end{split}
\end{align*}
For $p > \frac{1+n\Lambda}{1+2\Lambda}$ and $\Lambda \leq \frac{n-1}{n-2}$., take $\varepsilon_1$ small enough, then we obtain that
\begin{align*}
\begin{split}
\frac{2n-p-pn}{n-p}+\frac{(n-1)(p-1)}{(n-p)\Lambda}-p_{*}\varepsilon_1 < \frac{n(p-1-\frac{(2p-2)\Lambda-(p-1)}{\Lambda}+p\varepsilon_1)+p}{n-p}
\end{split}
\end{align*}
and
\begin{align*}
\begin{split}
0 \leq \frac{(2p-2)\Lambda-(p-1)}{\Lambda}+p\varepsilon_1 \leq p.
\end{split}
\end{align*}
If $\frac{p-1}{\Lambda}-2p+2-p\varepsilon_1  \leq \frac{2n-p-pn}{n-p}+\frac{(n-1)(p-1)}{(n-p)\Lambda}-p_{*}\varepsilon_1$, applying Corollary \ref{cor 3.12}, we obtain that
\begin{align}\label{4.1.1}
\begin{split}
\int_{\mathbb{R}^n} g^{-m}u^{2-p_{*}}W_{ij}W_{ji}\rho^{\gamma} &\lesssim R^{-2}\cdot R^{n-\frac{n-p}{p}(\frac{2n-p-pn}{n-p}+\frac{(n-1)(p-1)}{(n-p)\Lambda}) - \frac{n}{p}(\frac{(2p-2)\Lambda-(p-1)}{\Lambda})}\\
&= R^{\frac{p-1}{p\Lambda}-1},
\end{split}
\end{align}\[\]
taking $\varepsilon_1$ small enough and $R\rightarrow \infty$ in \eqref{4.1.1}, which yields $W_{ij}=0$. \\
If $\frac{2n-p-pn}{n-p}+\frac{(n-1)(p-1)}{(n-p)\Lambda}-p_{*}\varepsilon_1< \frac{p-1}{\Lambda}-2p+2-p\varepsilon_1$, applying Corollary \ref{cor 3.12}, we get that
\begin{align}\label{4.1.2}
\begin{split}
\int_{\mathbb{R}^n} g^{-m}u^{2-p_{*}}W_{ij}W_{ji}\rho^{\gamma} &\lesssim R^{-2}\cdot R^{n-\frac{n-p}{p-1}(\frac{2n-p-pn}{n-p}+\frac{(n-1)(p-1)}{(n-p)\Lambda}) - \frac{n-1}{p-1}(\frac{(2p-2)\Lambda-(p-1)}{\Lambda})}\\
&= R^{-2+\frac{3p\Lambda-n\Lambda-2\Lambda}{(p-1)\Lambda}},
\end{split}
\end{align}\[\]
taking $\varepsilon_1$ small enough and $R\rightarrow \infty$ in \eqref{4.1.2}. Since $\frac{3p\Lambda-n\Lambda-2\Lambda}{(p-1)\Lambda} <2$ for $p<n$, we also get $W_{ij}=0$.

(2). Case $p_n(\Lambda) < p < \frac{1+n\Lambda}{1+2\Lambda}$ and $\Lambda \leq 1+c(n)$ where $c(n)$ is a a number depending on $n$. In the following, we will present $c(n)$ respectively. Similar to the case 1, we only need to prove that
\begin{align*}
\begin{split}
\int_{\mathbb{R}^n} g^{-m}u^{2-p_{*}}W_{ij}W_{ji}\rho^{\gamma} \lesssim \int_{\mathbb{R}^n}g^{-m}u^{2-p_{*}}H^{2p-2}(\nabla u)\rho^{\gamma-2}|\nabla \rho|^{2}.
\end{split}
\end{align*}\[\]
Taking $m = \frac{p-1}{p\Lambda}-\varepsilon_1$ and observing that
\begin{align}
H(\nabla v) \leq (u^{p_{*}}g)^{\frac{1}{p}},
\end{align}
combining this with Remark \ref{rem 3.3} and Lemma \ref{lem 3.8}, then we get
\begin{align}\label{5.2.2}
\begin{split}
\int_{\mathbb{R}^n} g^{-m}u^{2-p_{*}}W_{ij}W_{ji}\rho^{\gamma}
&\lesssim \int_{\mathbb{R}^n}g^{-m}u^{2-p_{*}}H^{2p-2}(\nabla u)\rho^{\gamma-2}|\nabla \rho|^{2}\\
&\leq \int_{\mathbb{R}^n}g^{-\frac{p-1}{p\Lambda}+\varepsilon_1}u^{\frac{2n-p-np}{n-p}}(u^{p_{*}}g)^{\frac{2p-2}{p}}\rho^{\gamma-2}|\nabla \rho|^{2}\\
&= \int_{\mathbb{R}^n}g^{\frac{(2p-2)\Lambda-(p-1)}{p\Lambda}+\varepsilon_1}u^{\frac{np-3p+2}{n-p}}\rho^{\gamma-2}|\nabla \rho|^{2}\\
& \leq \int_{\mathbb{R}^n}g^{\frac{(2p-2)\Lambda-(p-1)}{p\Lambda}+\delta+\varepsilon_1}u^{\frac{np-3p+2-\delta p}{n-p}}\rho^{\gamma-2}|\nabla \rho|^{2}\\
& \lesssim \int_{\mathbb{R}^n}g^{\frac{p\Lambda-2\Lambda-(p-1)}{p\Lambda}+\delta+\varepsilon_1}u^{\frac{np-3p+2-\delta p}{n-p}}(u^{-p_{*}}H^{p}(\nabla u)+ u^{\frac{p}{n-p}})\rho^{\gamma-2}|\nabla \rho|^{2}\\
& \leq  R^{-2}\int_{\mathbb{R}^n}g^{\frac{p\Lambda-2\Lambda-(p-1)}{p\Lambda}+\delta+\varepsilon_1}u^{\frac{-2p+2-\delta p}{n-p}}H^{p}(\nabla u)\rho^{\gamma-2}
\\&+ R^{-2}\int_{\mathbb{R}^n}g^{\frac{p\Lambda-2\Lambda-(p-1)}{p\Lambda}+\delta+\varepsilon_1}u^{\frac{np-2p+2-\delta p}{n-p}}\rho^{\gamma-2}\\
& \leq R^{-2}\int_{\mathbb{R}^n}g^{\frac{p\Lambda-2\Lambda-(p-1)}{p\Lambda}+\delta+\varepsilon_1}u^{\frac{-3p+2-\delta p+n}{n-p}}H^{p-1}(\nabla u)\rho^{\gamma-3}|\nabla \rho|
\\&+ R^{-2}\int_{\mathbb{R}^n}g^{\frac{-2\Lambda-(p-1)}{p\Lambda}+\delta+\varepsilon_1}u^{\frac{-2p+2-\delta p+n-np}{n-p}}H^{p}(\nabla u)|W_{ij}|\rho^{\gamma-2}\\
&= I_1 + I_2.
\end{split}
\end{align}
For the first term $I_1$, we could take $\delta = \frac{n-3p+2}{p}+\varepsilon_1$. But for the convenience, we will use $\delta$ to calculate the following inequalities
\begin{align*}
\begin{split}
I_1 &= R^{-2}\int_{\mathbb{R}^n}g^{\frac{p\Lambda-2\Lambda-(p-1)}{p\Lambda}+\delta+\varepsilon_1}u^{\frac{-3p+2-\delta p+n}{n-p}}H^{p-1}(\nabla u)\rho^{\gamma-3}|\nabla \rho|\\
&\lesssim R^{2}\int_{\mathbb{R}^n}\left(u^{-p_{*}}H^{p}(\nabla u)+ u^{\frac{p}{n-p}}\right)^{\frac{p\Lambda-2\Lambda-(p-1)}{p\Lambda}+\delta+\varepsilon_1}u^{\frac{-3p+2-\delta p+n}{n-p}}H^{p-1}(\nabla u)\rho^{\gamma-3}|\nabla \rho|\\
&= R^{-2} \int_{\mathbb{R}^n}u^{-\frac{(n-1)(p\Lambda-2\Lambda-(p-1)+p\Lambda\delta)}{(n-p)\Lambda}-p_{*}\varepsilon_1 +\frac{-3p+2-\delta p+n}{n-p}}(H(\nabla u))^{\frac{2p\Lambda-3\Lambda-(p-1)+p\Lambda\delta}{\Lambda}+p\varepsilon_1}\rho^{\gamma-3}|\nabla \rho|\\
&+ R^{-2} \int_{\mathbb{R}^n}u^{\frac{p\Lambda-2\Lambda-(p-1)}{(n-p)\Lambda}+\frac{p}{n-p}\delta+\frac{p}{n-p}\varepsilon_1+\frac{-3p+2-\delta p+n}{n-p}}H^{p-1}(\nabla u)\rho^{\gamma-3}|\nabla \rho|\\
&= I_3 + I_4.
\end{split}
\end{align*}
For $I_3$ when $\varepsilon_1$ small enough, we observe that it satisfies the conditions of the Corollary \ref{cor 3.12},
\begin{align*}
\begin{split}
-\frac{(n-1)(p\Lambda-2\Lambda-(p-1)+p\Lambda\delta)}{(n-p)\Lambda} +\frac{-3p+2-\delta p+n}{n-p} \leq \frac{n(-\frac{p\Lambda-2\Lambda-(p-1)+p\Lambda\delta}{\Lambda}+p)}{n-p}
\end{split}
\end{align*}
and when $\Lambda \leq \frac{2n^2-11n+23+(n-3)\sqrt{4n^2-12n+57}}{4(n^2-5n+2)} = 1+c(n)$, it satisfies
\begin{align*}
\begin{split}
0 \leq \frac{2p\Lambda-3\Lambda-(p-1)+p\Lambda\delta}{\Lambda} \leq p.
\end{split}
\end{align*}
Thus applying the Corollary \ref{cor 3.12}, we obtain that
\begin{align*}
\begin{split}
I_3 &= R^{-2} \int_{\mathbb{R}^n}u^{-\frac{(n-1)(p\Lambda-2\Lambda-(p-1)+p\Lambda\delta)}{(n-p)\Lambda}-p_{*}\varepsilon_1 +\frac{-3p+2-\delta p+n}{n-p}}(H(\nabla u))^{\frac{2p\Lambda-3\Lambda-(p-1)+p\Lambda\delta}{\Lambda}+p\varepsilon_1}\rho^{\gamma-3}|\nabla \rho|\\
&\lesssim R^{-3}\cdot R^{n-\frac{n-p}{p}\left(-\frac{(n-1)(p\Lambda-2\Lambda-(p-1)+p\Lambda\delta)}{(n-p)\Lambda}-p_{*}\varepsilon_1 +\frac{-3p+2-\delta p+n}{n-p}\right)-\frac{n}{p}(\frac{2p\Lambda-3\Lambda-(p-1)+p\Lambda\delta}{\Lambda}+p\varepsilon_1)}\\
&= R^{\frac{p-1}{p\Lambda}-1},
\end{split}
\end{align*}
or
\begin{align*}
\begin{split}
I_3 &= R^{-2} \int_{\mathbb{R}^n}u^{-\frac{(n-1)(p\Lambda-2\Lambda-(p-1)+p\Lambda\delta)}{(n-p)\Lambda}-p_{*}\varepsilon_1 +\frac{-3p+2-\delta p+n}{n-p}}(H(\nabla u))^{\frac{2p\Lambda-3\Lambda-(p-1)+p\Lambda\delta}{\Lambda}+p\varepsilon_1}\rho^{\gamma-3}|\nabla \rho|\\
&\lesssim R^{-3}\cdot R^{n-\frac{n-p}{p-1}\left(-\frac{(n-1)(p\Lambda-2\Lambda-(p-1)+p\Lambda\delta)}{(n-p)\Lambda}-p_{*}\varepsilon_1 +\frac{-3p+2-\delta p+n}{n-p}\right)-\frac{n-1}{p-1}\left(\frac{2p\Lambda-3\Lambda-(p-1)+p\Lambda\delta}{\Lambda}+p\varepsilon_1\right)}\\
& = R^{-3}\cdot R^{\frac{-n\Lambda+4p\Lambda-3\Lambda+\delta p \Lambda}{(p-1)\Lambda}}\\
&= R^{-2}.
\end{split}
\end{align*}
Since $\Lambda = \lambda\max\limits_{\substack{\xi \in \mathbb{R}^n\\1 \leq i, j \leq n}}\left\{\frac{|\xi|^{2}(\nabla^{2}_{ij}H^{p}(\xi))} {p(p-1)H^{p}(\xi)}\right\}$, we obtain $I_3 \leq 0$ by $R \rightarrow +\infty$.\\
For $I_4$, we first take $\delta = \frac{n-3p+2}{p}+\varepsilon_1$ and observe that
\begin{align*}
\begin{split}
0 < \frac{-2p\Lambda+n\Lambda-p+1}{(n-p)\Lambda}< \frac{p\Lambda}{(n-p)\Lambda},
\end{split}
\end{align*}
when $\Lambda \geq 1$. Then we apply Corollary \ref{cor 3.12} to obtain that
\begin{align*}
\begin{split}
I_4 &= \int_{\mathbb{R}^n}u^{\frac{p\Lambda-2\Lambda-(p-1)}{(n-p)\Lambda}+\frac{p}{n-p}\delta+\frac{p}{n-p}\varepsilon_1+\frac{-3p+2-\delta p+n}{n-p}}H^{p-1}(\nabla u)\rho^{\gamma-3}|\rho_i|^{3}\\
&\lesssim  R^{-3}\cdot R^{n-\frac{n-p}{p}(\frac{p\Lambda-2\Lambda-(p-1)}{(n-p)\Lambda}+\frac{p}{n-p}\delta+\frac{p}{n-p}\varepsilon_1+\frac{-3p+2-\delta p+n}{n-p})-\frac{n}{p}(p-1)}\\
&= R^{-3+\frac{np\Lambda+2p\Lambda-n\Lambda+p-1-np+n}{p\Lambda}}.
\end{split}
\end{align*}
When $\Lambda \leq \frac{n^2-3n+3}{n^2-3n} = 1+\frac{3}{n^2-3n}$, we obtain $I_4 \leq 0$ by $R \rightarrow +\infty$.\\
For $I_2$, we use Young inequality yield that
\begin{align*}
\begin{split}
I_2 &= R^{-2} \int_{\mathbb{R}^n}g^{\frac{-2\Lambda-(p-1)}{p\Lambda}+\delta+\varepsilon_1}u^{\frac{-2p+2-\delta p+n-np}{n-p}}H^{p}(\nabla u)|W_{ij}|\rho^{\gamma-2}\\
&\lesssim \varepsilon_2 \int_{\mathbb{R}^n} g^{-\frac{(p-1)}{p\Lambda}+\varepsilon_1}u^{\frac{2n-p-np}{n-p}}W_{ij}W_{ji}\rho^{\gamma}\\
&+ \varepsilon_2^{-1}R^{-4} \int_{\mathbb{R}^n} g^{\frac{-4\Lambda-(p-1)}{p\Lambda}+2\delta+\varepsilon_1}u^{\frac{-3p+4-2\delta p-np}{n-p}}H^{2p}(\nabla u)\rho^{\gamma-4}\\
&= I_5 + I_6.
\end{split}
\end{align*}
We could choose $\varepsilon_{2}$ small enough and take $Q>1$ as a indicator of Young inequality. Then we only need to prove that
\begin{align*}
\begin{split}
I_2 &\lesssim \varepsilon_2^{-1}R^{-4}\int_{\mathbb{R}^n} g^{\frac{-4\Lambda-(p-1)}{p\Lambda}+2\delta+\varepsilon_1}u^{\frac{-3p+4-2\delta p-np}{n-p}}H^{2p}(\nabla u)\rho^{\gamma-4}\\
& \leq \varepsilon_2^{-1}R^{-4} \int_{\mathbb{R}^n} g^{\frac{-4\Lambda-(p-1)}{p\Lambda}+2\delta+\varepsilon_1}u^{\frac{-3p+4-2\delta p-np}{n-p}}\left(u^{p_{*}}g\right)^{2}\rho^{\gamma-4}\\
& = \varepsilon_2^{-1}R^{-4} \int_{\mathbb{R}^n} g^{\frac{2p\Lambda-4\Lambda-(p-1)}{p\Lambda}+2\delta+\varepsilon_1}u^{\frac{-5p+4-2\delta p+np}{n-p}}\rho^{\gamma-4}\\
& \leq \varepsilon_2 R^{-2} \int_{\mathbb{R}^n}g^{\frac{(2p-2)\Lambda-(p-1)}{p\Lambda}+\delta+\varepsilon_1}u^{\frac{np-3p+2-\delta p}{n-p}}\rho^{\gamma-2}|\rho_i|^{2}\\
&+ \varepsilon_2^{-2Q+1}R^{-2Q-2} \int_{\mathbb{R}^n}g^{\frac{(2p-2)\Lambda-(p-1)}{p\Lambda}+\delta+\varepsilon_1+(-\frac{2\Lambda}{p\Lambda}+\delta)Q} u^{\frac{np-3p+2-\delta p}{n-p}+\frac{-2p+2-\delta p}{n-p}Q}\rho^{-2Q+\gamma-2}\\
&= I_7+ I_8.
\end{split}
\end{align*}
Since $I_7$ is same as \eqref{5.2.2}, we only need to prove that $I_8 \leq 0$ when $R \rightarrow 0$.
\begin{align*}
\begin{split}
I_8 &= \varepsilon_2^{-2Q+1} R^{-2Q-2} \int_{\mathbb{R}^n}g^{\frac{(2p-2)\Lambda-(p-1)}{p\Lambda}+\delta+\varepsilon_1+(-\frac{2\Lambda}{p\Lambda}+\delta)Q} u^{\frac{np-3p+2-\delta p}{n-p}+\frac{-2p+2-\delta p}{n-p}Q}\rho^{-2Q+\gamma-2}\\
&= \varepsilon_2^{-2Q+1} R^{-2Q-2} \int_{\mathbb{R}^n}g^{\alpha} u^{\beta}\rho^{-2Q+\gamma-2},
\end{split}
\end{align*}
where $\alpha = \frac{(2p-2)\Lambda-(p-1)}{p\Lambda}+\delta+\varepsilon_1+(-\frac{2}{p}+\delta)Q$ and $\beta = \frac{np-3p+2-\delta p}{n-p}+\frac{-2p+2-\delta p}{n-p}Q$.
If we want to use Corollary \ref{cor 3.13}, we first check out that
\begin{enumerate}
  \item $0 < \alpha < 1 $
  \item $\frac{p}{n-p}\alpha + \beta< \frac{np-n+p}{n-p}$
  \item $-2Q-2+\max\{n- \alpha - \frac{n-p}{p}\beta, n-\frac{n-p}{p-1}\beta\}< 0$
\end{enumerate}
For item (1), direct calculation can be obtained that
\begin{align*}
\begin{split}
 \frac{p\Lambda-2\Lambda-(p-1)+\delta p\Lambda}{2\Lambda-p\Lambda\delta} < Q < \frac{(2p-2)\Lambda-(p-1)+\delta p\Lambda}{2\Lambda-p\Lambda\delta}
\end{split}
\end{align*}
For item (2),
\begin{align*}
\begin{split}
 Q > \frac{n\Lambda-2p\Lambda-(p-1)}{2p\Lambda}.
\end{split}
\end{align*}
For item (3),
\begin{align*}
\begin{split}
 Q < \frac{n-p-\delta p}{\delta p}.
\end{split}
\end{align*}
By observing that
\begin{align*}
\begin{split}
\frac{n\Lambda-2p\Lambda-(p-1)}{2p\Lambda} < \frac{p\Lambda-2\Lambda-(p-1)+\delta p\Lambda}{2\Lambda-p\Lambda\delta},
\end{split}
\end{align*}
since $\delta = \frac{n-3p+2}{p}$ and $p_n(\Lambda) < p < \frac{1+n\Lambda}{1+2\Lambda}$. Hence, we only need to calculate that the lower bound by
\begin{align}\label{5.2.3}
\begin{split}
\max\{\frac{p\Lambda-2\Lambda-(p-1)+\delta p\Lambda}{2\Lambda-p\Lambda\delta}    ,1\}    <  \frac{n-p-\delta p}{\delta p}.
\end{split}
\end{align}
A straightforward computation gives that \eqref{5.2.3} is equal to
\begin{align*}
\begin{split}
p > \frac{n+4}{5},
\end{split}
\end{align*}
and
\begin{align*}
\begin{split}
p > \frac{(5+n+3n\Lambda-2\Lambda)-\sqrt{(2\Lambda-n-5-3n\Lambda)^2-12(n^2\Lambda+n+2)}}{6}.
\end{split}
\end{align*}
Hence we can apply Corollary \ref{cor 3.13} and take $R\rightarrow \infty$ to get that
\begin{align*}
\begin{split}
I_8 \leq 0.
\end{split}
\end{align*}
When $W_{ij}=0$, this implies $Tr\{W^{2}\} = 0$ combining with Lemma \ref{lem 3.6}.  Hence $W = BF = TDQT^{-1} = 0$.
\vskip0.1cm
By the definition of $W_{ij}$, $W=0$ is equivalent to
\begin{align}\label{4.2.2}
\begin{split}
a^i = \lambda (x_i-(x_0)_i)u(x)^{p_{*}-1},
\end{split}
\end{align}
i.e.
\begin{align}\label{4.1.3}
\begin{split}
H^{p-1}(\nabla u)\nabla H(\nabla u)= \lambda (x-x_0)u(x)^{p_{*}-1},
\end{split}
\end{align}
which implies that
\begin{align}\label{4.1.4}
\begin{split}
 x-x_0 = \frac{1}{\lambda}u(x)^{1-p_{*}}H^{p-1}(\nabla u)\nabla H(\nabla u).
\end{split}
\end{align}
\vskip0.1cm
\noindent We notice that, acting $H_0$ on both sides of \eqref{4.1.4} and applying \eqref{2.1.7}, one could obtain that
\begin{align}\label{4.1.5}
\begin{split}
H_0(x-x_0) =  \frac{1}{\lambda}u(x)^{1-p_{*}}H^{p-1}(\nabla u)
\end{split}
\end{align}
\vskip0.1cm
\noindent Furthermore, according to \eqref{4.1.4} and \eqref{4.1.5} we have
\begin{align}\label{4.1.6}
\begin{split}
\nabla H(\nabla u) &= \frac{\lambda(x-x_0)u(x)^{p_{*}}-1}{H^{p-1}(\nabla u)}\\
&=\frac{x-x_0}{H_0(x-x_0)}.
\end{split}
\end{align}
\vskip0.1cm
\noindent Submitting \eqref{4.1.5} and \eqref{4.1.6} into \eqref{2.1.7}, and applying the property of 0-homogeneous of $\nabla H_0$ (the proof is same as $H_i$ in Lemma \eqref{lem 2.1}), then we compute
\vskip0.1cm
\begin{align}\label{4.2.15}
\begin{split}
\nabla u &= H(\nabla u)\nabla H_0\left(\nabla H(\nabla u)\right)\\
&=H(\nabla u)\nabla H_0(\frac{x-x_0}{H_0(x-x_0)})\\
&=\lambda^{\frac{1}{p-1}}u(x)^{\frac{p_{*}-1}{p-1}}H_0^{\frac{1}{p-1}}(x-x_0)\nabla H_0(x-x_0)\\
&=\frac{p-1}{p}\lambda^{\frac{1}{p-1}}u(x)^{\frac{p_{*}-1}{p-1}}\nabla \left(H_0^{\frac{p}{p-1}}(x-x_0)\right),
\end{split}
\end{align}
which implies that
\begin{align*}
\begin{split}
u = C_1 + C_2H_0^{-\frac{n-p}{p-1}}(x-x_0),
\end{split}
\end{align*}
\vskip0.1cm
\noindent for some $C_1, C_2 > 0$. Thus we have $u = U_{\lambda}$ and the proof of Theorem \ref{thm 1.1} is proved.
\end{proof}

\medskip
\vskip0.1cm
\noindent\textbf{Open problem:}
\medskip
Prove that for $n \geq 2$ and for $1 < p< n$, assume that $u \in W^{1,p}_{loc}(\mathbb{R}^n)$ is a positive weak solution of \eqref{1.5}. Then $u$ must take the form as
\begin{equation*}
U_{\lambda}(x) = \left(\frac{(\lambda^{\frac{1}{p-1}}(n^{\frac{1}{p}}(\frac{n-p}{p-1})^{\frac{p-1}{p}})}{\lambda^{\frac{p}{p-1}}+H_0(x)^{\frac{p}{p-1}}}\right)^{\frac{n-p}{p}},
\end{equation*}
up to some translation.

{\bf Data availability:} Data sharing not applicable to this article as no datasets were generated
or analysed during the current study.

\end{document}